\documentclass[12pt,dvips]{amsart}
 
\usepackage[dvips]{graphics}
\usepackage{epsfig}
\date{October, 1997}
 
\newcommand{\A}{{\mathbb A}}
\newcommand{\C}{{\mathbb C}}
\newcommand{\D}{{\mathbb D}}

\newcommand{\F}{{\mathbb F}}

\newcommand{\PP}{{\mathbb P}}
\newcommand{\Q}{{\mathbb Q}}

\newcommand{\Z}{{\mathbb Z}}
\newcommand{\uppermin}[1]{#1^{\text{min}}}
\newcommand{\tablehead}{
 \begin{tabular}{|p{0.2in}|p{1.2in}|p{0.8in}|p{0.5in}|p{1.4in}|p{1.2in}|} 
 \hline
 &
 \centering {configuration} &
 \centering {$B_1$} &
 \centering {max($b$)} &
 \centering {complements} &
 \begin{minipage}{1.2in}
 \begin{center}  {divisors\\(Pic(S)=$\Z[H]$)} \end{center} \end{minipage}\\  

 \hline 
}
\newcommand{\settablewidth}{
 \begin{tabular}{|p{0.2in}|p{1.2in}|p{0.8in}|p{0.5in}|p{1.4in}|p{1.2in}|} 
 \hline
} 
\newcommand{\includefigure}[1]{
 \begin{minipage}{1.2in}\begin{center}
  \vspace{0.1in}
  \includegraphics{elliptic.ex#1.eps}
  \vspace{0.1in}
 \end{center}\end{minipage} 
}
\title{Classification of Exceptional Complements:\\Elliptic Curve Case}

\author{Terutake Abe}
 
\address{\centering Department of Mathematics\\
The Johns Hopkins University\\
3400 N. Charles Street \\
Baltimore, MD., 21218, U.S.A.\\}
 
\email{abe@chow.mat.jhu.edu}
 
\thanks{I would like to thank Professor Shokurov for setting the problem
and for his valuable suggestions.}

\newtheorem{theorem}{{\sc Theorem}}[section]

\newtheorem{lemma}[theorem]{{\sc Lemma}}

\newtheorem{case}{Case}

\newenvironment{claim}{\smallskip\noindent{\bf Claim:}}{\smallskip}
\renewcommand{\subsection}[1]{\bigskip\noindent{\bf #1}\\}  
\begin{document}
 
\begin{abstract} 
We classify the log del Pezzo surface $(S,B)$ of rank 1 with no  
1-,2-,3-,4-, or 6-complements with the additional condition that
$B$ has one irreducible component $C$ 
which is an elliptic curve and $C$ has  
the coefficient $b$ in $B$ with $\tfrac1n\lfloor(n+1)b\rfloor
=1$ for n=1,2,3,4, and 6.
\end{abstract}

\maketitle
%
%
%
%

\section {Introduction}

This paper is a part of the project to classify 
``log del Pezzo surfaces with no regular complements'', that is, 
the pairs (S,B) of surface S and boundary B on S such that:
\newcounter{ex}
\begin{list}{(EX\arabic{ex})}{\usecounter{ex}}
\item  $-(K+B)$ is nef ($(S,B)$ is ``quasi log del Pezzo''), 
\item  $-(K+B)$ has no regular complements i.e. it has no n-complements 
for any of $n\in \{1,2,3,4,6\}$.
\end{list}
\noindent
We assume throughout that coefficients of $B$ are ``standard'', i.e.
$B= \sum b_iC_i$ with $b_i = \frac{m-1}m$ where $m$ natural number, or
$b_i \ge \frac67$. 
An invariant $\delta$ for such a pair is defined in [Sh2, 5] by
\begin{align*}
\delta(S,B) = \sharp\{E| E  
                 & \mbox{ is an exceptional or non-exceptional divisor} \\
                 & \mbox{ with log descrepancy } 
                        a(E) \le \tfrac17 \mbox{ for } K+B\}
\end{align*}
and it was proved there that $\delta \le 2$ ([Sh2, Th.5.1]).
We can assume, after crepant blow ups of
exceptional $E$'s with $a(E)\le \frac17$, that  those $E$ are all 
non-exceptional, and thus,
\begin{list}{(EX3)}{}
\item   $(S,B)$ is $\frac17-$ log terminal.
\end{list}
Now define the divisor $D$ by 
$D= \sum d_iC_i$ where  $d_i=1$ if $b_i\ge \frac67$ and $d_i=b_i$ otherwise.
And write $C = \lfloor D\rfloor = 
\sum\limits_{a(C_i)\le\frac17}C_i$. We know by [Sh2, Lemma 4.2] that
if $\delta \ge 1$ then we can successively contract
curves semi-negative  with respect to $K+B$, but not components of $C$, and 
thereby assume 
\begin{list}{(EX4)}{\usecounter{ex}}
\item  $\rho(S)=1$.
\end{list}
The conditions (EX1),(EX2) and (EX3), as well as the condition on the 
coefficients, are preserved under this reduction.
We form a minimal resolution $f:(S^{\text{min}},B^{\text{min}}) \to (S,B)$ 
where $B^{\text{min}}$ is a crepant pullback, i.e.
$K_{S^{\text{min}}}+B^{\text{min}}= f^*(K+B)=K+f^{-1}(B)+\Sigma e_jE_j$ 
satisfies 
$K + B^{\text{min}} \cdot E_j = 0$ for all $j$. 
From $S^{\text{min}}$ we contract $(-1)$-curves successively to 
get a smooth model $S^\prime$ which is either $\PP^2$ or $\F_m$:
$$ g:(S^{\text{min}},B^{\text{min}})\to(S^\prime,B^\prime).$$
If $\delta \ge 1$ we have $p_a(C) \le 1$, 
and the same is true for the birational
image of $C$ on $S^\prime$ as well([Sh2, Prop.5.4]).
In this paper we consider the case 
\begin{itemize}
\item $\delta = 1$. 
\end{itemize}
Thus, $C$ is an irreducible curve of arithmetic genus $\le$ 1.
We write $C=C_1$, and $B=bC+\sum_{i=2}^r b_iC_i = bC+B_1$, 
with $b_i\in\{0,\frac12,\frac23,\frac34,\frac45,\frac56\}$ and 
$b\geq\frac67$.

%
%
    \subsection{The dual graph of a configuration}
%
%

In the following we use the language of graphs to
talk about the the configuration of curves.

The dual graph of  a configuration of curves is a (weighted-multi) graph 
where we have a vertex
for each curve and an n-ple edge for each
intersection point with multiplicity $n$ between two curves.
Each vertex has a weight 
$\in \Z$ which is the self-intersection number 
of the curve.

Graphically, we use $\bullet$ (``b(lack)-vertex'') to represent exceptional
curves with self-intersection number $\le -2$, $\circ$ 
(``w(hite)-vertex'') for
$(-1)$-curves, and squares for curves
with non-negative self-intersection. The weight of a vertex is shown
by a number next to each vertex, and multiplicity of an edge by
the number of lines joining the two end vertices.

``Blow up of an edge'' means the transformation of the 
graph reflecting the blow up of the corresponding point, that is,
introduce a new white vertex, decrease the multiplicity of edge
by 1, decrease the weight of the both 
end  vertices of the  edge by 1, and join them to the new white vertex by
a simple edge.
``Blow up of a vertex'' reflects the blow up of a point on
the curve outside the intersection with neighboring curves:
introduce a white vertex, decrease the weight of the vertex by 1,
and join it to the new white vertex by a simple edge.
Blow up of a complete subgraph of any cardinality $k$ can be defined
in the same way.

%
%
               \subsection{Types of Singularities on $C$}
%
%

\begin{lemma} 
\begin{itemize} 
\item[(i)]  The singularity of $C$ is at worst a node, and it is 
          outside $Sing(S) \cup Supp(B_1)$.
\item[(ii)] At most one component $C_i$ of $B_1$
          passes through each point $P\in C$. If $P$ in a smooth point of
          $S$, then the intersection is normal 
          ,with one possible 
          exception where $C_i$ has coefficient 
          $\frac12$ and has a simple tangency with $C$ at a smooth point 
          $P$ of $S$.
\item[(iii)] Singularity $P$ of $S$ on $C$ is a cyclic quotient  
          singularity, i.e. log terminal singularities with resolution graph
          $\A_n$ (a chain), where  $C$ meets one end curve of the chain
          normally.
          If another component 
          $C_i$ passes through $P$, then it 
          meets the other end curve normally. 
\end{itemize}
\end{lemma}

\proof

Note that $K+D$, as defined above, is log canonical by the existence 
of local complements
([Sh1, Cor.5.9.]).  Then all the statements follow from the 
classification of surface log canonical
singularities ([Ka],or [Al2]) and $\frac17$-log terminal condition.   
For example, for (iii), if we had a
type $\D_n$ singularity, (case (6) in [Ka, Th.9.6]) we would have 
a log decrepancy $\le \frac17$. Note also that the exception 
in (ii) is the only case where $K+D$ is not log terminal at $P$ 
([Sh2, Prop.5.2]).   

\qed

As is well known, the singularities mentioned above are isomorphic,
analytically,  to the origin 0 in the quotient of 
$\C^2$ by the action of cyclic group $\mu_m$ of order $m$, 
where the generator 
$\varepsilon = e^{\frac{2\pi i}{m}}$ acts by
$$ (z_1, z_2) \mapsto (\varepsilon^{-k}\cdot{z_1}, \varepsilon\cdot z_2),
\text{   where  } 1 \le k \le m \text{    and    }  gcd(m,k)=1.$$
The minimal resolution of such a singularity has 
a chain of rational curves
$E_1,E_2,\cdots,E_r)$ as its 
exceptional locus, and the coe
continued fraction expansion
$$\frac{m}{m-k} = w_1 - \cfrac{1}{w_2-
       \cfrac{1}{w_3-\dotsb}}$$

\noindent
give their self intersection numbers ({\it cf.} for example, [Ful]).
We call such a singularity $P$ type $[m,k]$. We extend this correspondence 
to incorporate the information on the
component $C_i$ that passes through $P$ ({\it cf.} [Sh1,Cor.3.10],
[Sh2,Lemma 2.22]).
 
Namely, if the component $C_i$ has the standard
coefficient $\frac{d-1}{d}$ and the singularity $P$ has type $(m^\prime,
k^\prime)$,  
we represent it by the pair $(m,k)=(dm^\prime,dk^\prime)$.
The ``dual graph'' of the minimal resolution of this singularity is:

\begin{figure}[h] \begin{center}
\includegraphics{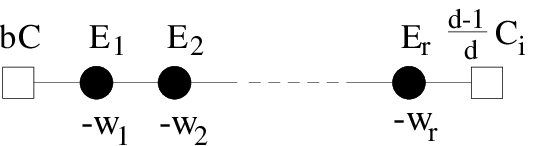}
\caption {}
\end{center} \end{figure}
\noindent

Generalizing the notation of [KM], we may denote the same singularity
by $(\underline{w_1},w_2, \cdots, w_r)_d$
with the underline indecating the curve meeting $C$. 

This singularity has the minimal log discrepancy
$$mld(P,K+B)=a(E_1)=\tfrac{1+(m-k)(1-b)}{m} \text{ }(\le \tfrac{1+\frac17(m-k)}{m}), $$
where $m = d \cdot \text{(index of $P$)}$. Also we denote the co-discrepancy,
or the coefficient, of $P$ by $e(P,K+B)=1-mld(P,K+B)$. 

Now the $\frac17$-log terminal condition
$$\tfrac{1+\frac17(m-k)}{m} > \tfrac17$$
is equivalent to $k < 7$. Therefore the possible singularities on $C$ 
are put into 21 $=6(6+1)/2$ (infinite) series according to the pair 
$(m (\text{mod }k),k)$ with $1\le k\le 6$. 
This will be convenient later on.

%
%
%
                 \section {Elliptic Curve Case}
%
%
%

Now we start the classification of the case $p_a(C) = 1$. Thus,
$C\in S$ is a smooth curve of genus 1 or a rational curve with one node.
We call it the ``elliptic curve case''.
\begin{lemma}
In the ``elliptic curve case'', 
the condition (EX2) is equivalent to the condition that $(S,B)$ 
has log-singularities on $C$. That is,
either $S$ has singularities  on $C$, or $B$ has components other than
$C$ (which intersect $C$ since $\rho(S)=1$).
\end{lemma}
\proof
If $(S,B)$ is smooth on $C$ , then $(K+f^*(D)).C=(K+C).C 
= 0$ on $\uppermin{S}$,
so $K+D=K+C\sim 0$ on $S$ and (EX2) is not satisfied. In fact
$K+B=K+bC$  has a 1-complement.
On the other hand if $(S,B)$ has a singularity on $C$, 
then $(K+f^*(D)).C>(K+C).C=0$ so we have $K+D>0$ on $S$, 
which implies (EX2).
\qed
 
The case when $S$ is a cone ($\PP^2 \text{ or } \Q_m$) has been classified
elsewhere and from it we have only one case with $C$=elliptic: $S=\F_2,
C=$double section,$B_1=\frac12C_2$ where $C_2$ is a generator of the cone.
Then $C\equiv 2H \equiv -K, C_2\equiv \frac12H$. So
$K+\frac67C+\frac47C_2\equiv 0$ and $K+B$ has 7-complement $= 0$. 
It also has the trivial 8-complement: $K+\frac78 C+ \frac12 C_2 \equiv 0$
(This is the entry $\sharp$1 in the table at the end).

From now on we assume $S$ is not a cone.

\begin{lemma} $C^2 \ge 3$ on $S^{\text{\rm min}}$. If $(S,B)$ has 
two singularities on $C$, then $C^2 \ge 6$.  
On the other hand,
the minimum log discrepancy of the singularity $P$ on $C$ with respect 
to $K+B$ (hence also with respect to $K+bC$) is at least $1-(C^2/7)$.
\end{lemma}
\proof

Because $-(K+B)$ is nef,
\begin{align*}
 0\geq (K+B).C &= K_{S^{\text{min}}} +bC +\Sigma b_iC_i +
                  \Sigma d_jE_j \cdot C\\
           &\geq  -(1-b)C^2 + \Sigma b_i + \Sigma_{P}(1-mld(P))\\ 
           &\geq  -(1-b)C^2 + min\{b_i, 1-mld(P)\} \\
           &\geq  -\tfrac17 C^2 + \tfrac37. 
\end{align*}
Note that, because of Lemma 1.1, $1-mld(P)=d_j$ 
for the exceptional curve $E_j$ meeting $C$.
The last inequality holds because  
we have at least one nonzero $b_i$ or $1-mld(P)$  by Lemma 2.1 
and the minimum nonzero value for $b_i$ is $\frac12$, that for $1-mld(P)$ is 
$\frac12 \cdot \frac67 = \frac37$, 
the latter being attained when $P$ in duVal of type $\A_1$.  
Therefore, $C^2 \geq 3$.  By the same calculation, 
if there are two 
singularities on $C$ we have $0\geq -\frac17 C^2 + \frac67$. 
On the other hand, the second inequality in particular implies that
$1-mld(P) \le (1-b)C^2 \le \tfrac17 C^2$, whence the second assertion. 
$\qed$

%
%
%
                \subsection{Reduction to $\F_2$}
%
%
%

We need the following
\begin{lemma} Let $E$ be a $(-1)$-curve  on $\uppermin{S}$. Then  
on its image $f_*(E)$,
$S$ has either at least two singularities, or one singularity
that is not log-terminal for $K+E$.  
\end{lemma} 
\proof If , on $E$, $S$ had at most one singularity $P$ that is 
log-terminal for $K+E$, i.e. a cyclic quotient singularity
such that $E$ meets one end curve $E_1$ of the chain of the resolution,
then we would have
$$ (f_*(E))^2=E.f^*f_*(E)=E.(E+(1-a(E_1))E_1)=-1 + (1 -a(E_1)) 
< 0$$
which is absurd since $\rho(S)=1$. 
\qed

\medskip
Now we can prove the
\begin{lemma} We can always obtain $\F_2$ as a smooth model of $S$ (and 
$C$ as a double section).
\end{lemma}
\proof

$p_a(C)=1$ means that after reconstruction, $C$ is either a
cubic in $\PP^2$, curve of bidegree $(2,2)$ on $\F_0$, or a double section
of $\F_2$. Suppose $S^\prime$ is $\PP^2$ and $C$ is a cubic, since there
are no irreducible curves with arithmetic genus 1 on $\F_m$, with $m\ge 3$.
If $g: S^{\text{min}} \to \PP^2$
contracts two or more exceptioncal curves to a point $P \in \PP^2$, 
then we can choose different contractions to get $S^\prime = \F_2$.
Therefore we may assume that we have  only one exceptional curve for $g$ over
each center $P \in \PP^2$, and we shall derive a contradiction.

Since all the curves contracted by $g$ are $(-1)$-curves on 
$S^{\text{min}}$, no exceptional curve $E_i$ for $f$
are contracted and all of them are present on $\PP^2$ as divisors. 
Thus we have an inequality
$$ 
0 \geq deg(K_{\PP^2}+bC+B_1^\prime) 
  = -3+\tfrac67\cdot 3 + deg(B_1^\prime)
  = -\tfrac37 + deg(B_1^\prime) $$

Therefore, since the coefficients are standard, no component $C_i$ 
other than $C$ are present on $\PP^2$. And we have
\begin{equation} \Sigma d_j \le \Sigma d_j\cdot deg(E_j) = deg(B_1^\prime)
\le \tfrac37 
\tag{$*$} \end{equation}  

We have the two possibilities:

(1) $B$ has at least one  component, say $C_2$, other than $C$. Then by 
the above, $C_2$ must be contracted on $\PP^2$ and is a 
$(-1)$-curve on $S^{\text{min}}$.  Therefore, by Lemma 2.3,
$S$  must have
either at least two singularities on $C_2$, or a singularity that is
not log-terminal for $K+C_2$.  
In the former case, then,
we would have $\Sigma d_j \ge (\frac12 + \frac12)b_2 \ge \frac12 > \frac37$,
contradicting ($*$).
In the latter case, we have an exceptional curve $E$ with 
$a(E,K+C_2)\le 0$. Then because $a(E, K+b_2C_2)$  is a linear function 
of $b_2$ and we also have $a(E,K+0\cdot C_2)=a(E,K) \le 1$, 
we have   $a(E,K+b_2 C)\le 1-b_2$.
Thus $d_2 = 1-a(E,K+B) \ge 1-a(E,K+b_2C_2) \ge b_2 \ge \frac12 > \frac37$,
again a contradiction to ($*$).

(2) $B$ has no other components than $C$, i.e. $B=bC$, 
and $S$ has a singularity on $C$. Then ($*$) implies 
that 
and we have $K+B > 0$ except in the following case: $S$ has only one  
duVal singularity $P$ of type $\A_1$ on $C$, 
the exceptional curve $E_1$ of the resolution of $P$ is a line on $\PP^2$,
$b=\frac67$, and $B_1^\prime$ has no other component than $E_1$, so that  
$K_{\PP^2}+ B^\prime = K + \frac67 C + \frac37 E_1 \sim 0$. 
In particular all 
the singularities on $S$ are duVal so 
$$f^*(K+B)=K_{\uppermin{S}}+\uppermin{B} 
= K+\frac67 C+ \frac37 E_1$$
Also, the triviality of $K+B$ means that pull back $g^*$ is crepant
so that the above is also equal to
$g^*(K+B^\prime).$
On the other hand,
since $E_1.C = 3$ on $\PP^2$ and $E_1.C=1$ on $S^{\text{min}}$,
two of the intersection points of 
$E_1$ and $C$ has to be blown up on $S^{\text{min}}$. The exceptional
curve $E$ for the first of such  blowups would have the coefficient 
$\frac67 + \frac37 - 1 = \frac27$ 
in $K_{S^{\text{min}}}+B^{\text{min}} = g^*(K+B^\prime)$. 
Contradicting the explicit form of $\uppermin{B}$ given above. 

If we have a model $S^\prime = \F_0$, then we have had 
at least one contraction
of $(-1)$-curve so we can get $S^\prime = \PP^2$ by choosing other 
contractions, and we are reduced to the previous case.

\qed

Therefore we have a $\PP^1$-fibration $p:S^{\text{min}}\to \F_2\to \PP^1$.
Now our strategy for the classification is to start from $\F_2 = S^\prime$,
make blow ups to construct $S^{\text{min}}$, choose $B^{\text{min}}$
on it so that resulting $(S,B)$ would have singularities on $C$ 
($\Leftrightarrow$ (EX2) by Lemma 2.1) and 
would satisfy (EX1),(EX3), and (EX4). 

The conditions (EX1) and (EX4) implies that  the number of $(-n)$-curves,
$n\ge 2$, on $S^{\text{min}}$ must equal $\rho(S^{\text{min}})-1$, and they
are all exceptional for the resolution $f$. 
These curves are either in the fibres of $p$, or they are not, i.e.
they are (multi-) sections of $p$.  As for the number of curves of each type,
we have the following:

\begin{lemma} {\rm (\lbrack Zhang, Lemma 1.5 \rbrack)} We have
\begin{eqnarray*}
r & = & \sharp\text{\{Exceptional curves $E_i$'s of the resolution f 
        that are not in the fibres of $p$\}} - \text{ 1 } \\
  & = & \sharp\text{\{$(-1)$-curves on $S^{\text{min}}$ that are 
        in the fibres of $p$.\}} \\
  & & -\ \  \sharp\text{\{Singular fibres of p\}}
\end{eqnarray*}
\end{lemma}
\proof

Add ($2 + \sharp\{E_i$'s that are  in the fibres of $p$\}) to both sides,
and we get two expressions for $\rho(S^{\text{min}})$.$\qed$\

%
%
%
        \subsection{The search for exceptions}
%
%
%

\begin{case}
$r=0$, i.e. minimal section $\Sigma$ on $\F_2$ is the only $E_i$ 
with $p(E_i) = \PP^1$.
\end{case}
Then there is only one $(-1)$-curve in each singular fibre of $p$. Therefore
on each fibre $F$ modified we have to have initially two blow ups at the same 
point $P$. Suppose $C^2 = w$ before the modification, then
according as the intersection multiplicity $i=I(P;F\cap C)=2,1,$ 
or $0$, i.e. according  as $P =$ tangency of $F$ and $C$, 
normal intersection of $F$ and $C$, 
or $P \in F \setminus (F\cap C)$, 
we get one of the three dual graphs in the Figure 2 below.

\begin{figure}[h]
\begin{center}

\begin{tabular}{p{1.0in}p{0.5in}p{1.0in}p{0.5in}p{1.0in}} 
\includegraphics{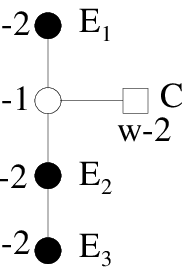} & \hspace{0.5in} &
\includegraphics{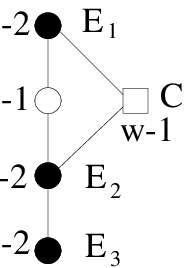} & \hspace{0.5in} &
\includegraphics{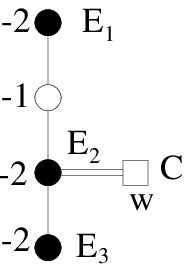}  \\
\mbox{(I)} &  & \mbox{(II)} & & \mbox{(III)} 
\end{tabular}
\caption{}
\end{center}
\end{figure}

In the figure the b-vertex at the bottom is the minimal
section $\Sigma \in \F_2$. 
In the case (III),  the curve C
and neighboring $(-2)$-curve ($=F$) have either two normal intersections,
one simple tangency, or $C$ has a node on $F$.  

Case(III) gives a non log canonical point ({\it cf.} Lemma 1.1(i)) 
and is excluded.
Case(II) gives one example with trivial complement 
(entry $\sharp2$ in the table at the end):
\begin{align*}
S &= \text{Gorenstein del Pezzo surface with singularities $\A_1+\A_2$}, \\
C &= \text{elliptic curve through $\A_1$ and $\A_2$ points},\\
K+&B = K+\frac67C\equiv 0 \\
7(K_{S^{\text{min}}}+ & B^{\text{min}})= 7(K + \frac67C + 
\frac37E_1 + \frac47E_2 + \frac27E_2) \sim 0
\end{align*}

\noindent
(Following [MZ],  we denote the Gorenstein del Pezzo surfaces of rank 1
by its singularity type, for example, $S(A_1+A_2)$ for the surface above,
and their resolution by e.g. $\tilde{S}(A_1 + A_2)$.)

Since we already have $K+ B \equiv 0$, if we make any more blow ups
(which have to be on the unique $(-1)$-curve) or
add other components to $B$, we would have $K+B > 0$ and $(S,B)$ will 
violate (EX1). So we need not 
consider this case any longer.
Thus we are left with case (I), i.e. two initial blow ups at
the ramification point of $C \to \PP^1$ (tangency of $C$ and a fibre). 
In particular, in all the remaining cases, $C^2 \le 6$, because $C^2=8$ 
on  $\F_2$.

This implies that a smooth fibre $F$ cannot be a component of $B_1^\prime$, 
because if it were, we would have
$0 \ge (K+bC+B_1).C \ge -(1-b)C^2 + \frac12 F.C 
\ge -\frac67+\frac12\times 2 = \frac17$,
a contradiction.  Therefore only singularities on $C$ 
are those coming from the intersection
of $C$ and the singular fibres.  

After (I), we can only  blow up 
a point on the unique $(-1)$-curve on each fibre: otherwise we would introduce 
more than one $(-1)$-curves in a fibre, violating $r=0$. 
There are two types of such blow ups. One is the blow ups 
of the intersection of $C$
and the $(-1)$-curve, (blow up of the edge between the white vertex and $C$)
which  decrease $C^2$. 
The other is the blow ups of a point of $(-1)$-curve  outside $C$.

We start from  the first type of blow ups and get the resolutions of 
Gorenstein log del Pezzos of rank 1 with $K^2 = C^2 \ge 3$ (Lemma 2.2):
$$\begin{array}{lllllll}
\tilde{S}(A_1 + A_2)&\longrightarrow&
\tilde{S}(A_4) &\longrightarrow&
\tilde{S}(D_5) &\longrightarrow&
\tilde{S}(E_6)\\
 &\searrow&                     &\searrow&                   & &\\
 &        &\tilde{S}(2A_1+A_3)  &\longrightarrow&
\tilde{S}(A_1+A_5).& &
\end{array}$$

Each ``$\longrightarrow$'' represents one blow up, and each ``$\searrow$''
two blow ups on a new fibre .

Then, starting from one of these, we make the second type of blow ups,
which decrease the minimal log discrepancy of $S$, 
until either (EX1) or (EX3) is violated (see below).  
The Gorenstein rank 1 surfaces listed
above are the image of $\uppermin{S}$ under the morphism $\phi_{|C|}$ 
defined by the linear system $|C|$ on it. 
We denote it by $S_C$, and its resolution
(one of the above) by $\tilde{S_C}$. 
 
Note that $C$ meets every $(-1)$-curve $E$ on $S_C$ since
$C\sim -K_{S_C}$ and $-K.E=1$.  Consider blow ups on one fibre 
starting at one such $E$. 
By Lemma 2.3, on $E$, $S$ has either at least two 
singularity or one singularity that is not log-terminal for $K+E$.
That is, on $\tilde{S_C}$, either $E$ meets at least two trees $T_1,T_2$
of b-vertices,
or one tree $T_3$ that gives non-log-terminal point for $K+E$. 

Now consider the transformation of the subgraph consisting or
$C$, $E$, and trees of b-vertices $T_i$ meeting $E$ on $S_C$.
It should always contain a unique w-vertex.

If we blow up the vertex $E$, i.e. blow up 
a point on $E$ other than the intersection points
with neighboring exceptional curves, then after the transformation 
$C$ would meet the b-vertex $E$ in the black graph $T_1 - E - T_2$ or 
$E - T_3$. 
Either of these would contracts to a
non-log-terminal point on $S$ for $K+C$, contradicting Lemma 1.1.   
(For an example  of the first situation, consider blow up of the
white vertex in the configuration (I) in the Figure 2 above.
For the second, consider the same in the configuration of the table
$\sharp 9$.)  Therefore the first blow up has to be at the intersection
point of $E$ and one of the neighboring b-vertices, i.e. blow up
of the edge joining $E$ and one of its neighbors. 

The  same argument, repeated for the new white vertex $E_1$ at each stage, 
shows that successive blow ups also must be  
at the edge joining $E_1$ and a neighboring b-vertex, because the 
trees now meeting $E_1$ are even bigger than those that met $E$.
Thus, by induction, we see that the full inverse image of 
$E$ is of the form $E - T - E_1 - T^\prime$, 
where $T$ and $T^\prime$ are chains of b-vertices ($T$, or $T^\prime$
may be a part of a larger tree. And $E$ may meet another tree
$T^{\prime\prime}$ in which case $T$ should be empty --- Remember that
$C$ meets $E$), and $E_1$ is a w-vertex. 
The blow up described above either  increases the weight of an end vertex
of $T$ next to $E_1$, 
or adds one $(-2)$-curve $E_1$ to it, 
depending on which side of $E_1$ we blow up.
Either of such  tranformations
(those which  preserve log-terminal property), 
if repeated infinitely many times,
make the log-discrepancy with respect to $K+bC$ 
of the resulting singularity on $C$ monotonically 
decrease toward $1-b \le \tfrac17$.  Hence
by Lemma 2.2, after finite number of steps, (EX3)
will be violated (or perhaps, (EX1) may be violated first).  
Therefore this procedure of successive blowups 
must terminate. 

We can now refine the lemma 2.2 as follows:

If $C^2 < 6$ then we have only one singularity by Lemma 2.2. 
But on the other hand, if $C^2 \ge 5$ we can have only one singular fibre, 
which means that in every case we have only one singularity of 
$(S,B)$ on $C$.  (EX1) restricts the possible
types of singularities $[m,k]$ on $C$ as follows:
\begin{eqnarray*}
0\geq K+bC+B^\prime \cdot C &=& -(1-b)C^2+\tfrac{(k-1)+b(m-k)}{m}\\
                            &=& \frac17C^2+\tfrac{(k-1)+\frac67(m-k)}{m},\\
                            &or&\\
                   (6-C^2)m&\leq& 7-k.
\end{eqnarray*}
  
In this way, we find that there 
are 20 possible $S$'s , with a few different $B$'s for some of the $S$'s. 
These are summarized in
the table below.

\begin{case} $r=1$, i.e. we have one exceptional curve, say $E$, other than
$\Sigma$ that is a section of $\PP^1$-fibration $p$.
\end{case}
Thus, exactly one fibre contains two $(-1)$-
curves in it.  If we modify at any other  fibre it has to
start like (I) of the Figure 2 (two blow ups at the tangency with the fibre)
because (II) and (III) have
been elimineted. In particular each time we blow up on a new fibre
we decrease $C^2$ by at least 2. 

\begin{claim} Any exceptional curve $E$ that is a (mutli-)section 
of $p$ is in fact a 1-section that is disjoint from $\Sigma$.
\end{claim} 
\proof
Let $E$ be a (multi-)section, and $d=\text{mult}_E(B_1^\prime)$.
Then if $F$ is a fibre of $p$, we have
$$0 \ge (K+B^\prime).F \ge (K+bC+dE).F \le -2+\tfrac67\cdot 2 + d.$$
Hence $d \le \frac27 < \frac37$.  
So $E$ cannot intersect $C$ on $\uppermin{S}$. Therefore all the intersection
point of $C$ and $E$ have to be blown up on $\uppermin{S}$.
If $E$ is not a 1-section disjoint from $\Sigma$, then we have $C.E\ge 6$ 
on $\F_2$. So we would have $C^2 \le 8 - 6 =2$ on $\uppermin{S}$ 
which is impossible according to the lemma 2.2.  This proves the claim.

So let $E$ be a simple section disjoint from $\Sigma$. Then $E.C=4$.

Suppose $E$  intersects $C$ at one point
with multiplicity 4. Then after four blowups at this point 
we get $\tilde{S}(A_1+A_3)$ (the configuration of the table $\sharp 17$,
with a different choice of fibration), 
which has already been studied in the case 1 above.

If $E$ intersects $C$ at two points with multuplicity
3 and 1  respectively, then by the above observation 
we have at least $3+2=5$ blow ups on $C$, which gives 
$\tilde{S}(3A_2)$, with $C$ passing through three
$(-1)$-curves joining three $\A_2$ points. Since $C^2=3$
by Lemma 2.2 $C$ can have at worst
$\A_1$ (=``type [2,1]'') point on it, 
but that cannot be attained: We could at best choose $B_1=\frac12 C_2$ where 
$C_2=$ (image of one of the $(-1)$-curves meeting $C$) and thus get 
type [2,2] point on $C$, which
is worse than $\A_1 = type [2,1]$. 

If $E$ intersects $C$ at more than 3 points,
then we have at least six blow ups on $C$ thus  $C^2 \le 2$, 
which is impossible by Lemma 2.2.

Finally, if $E$ intersects $C$ at two points with multiplicity
2 each, we would have two $(-1)$-curves in each fibre,
and this violates (EX4). 

Thus, we get no new examples from case 2. 

\begin{case} $r\ge 2$, i.e. we have at least two exceptional
curves, say, $E_1$ and $E_2$, other than $\Sigma$, that are sections of $p$.
\end{case}
$E_i$ are simple sections. Then because $C.E_i=4$ and 
$E_1.E_2=2$, we must have at least $4+4-2=6$ blow ups on $C$ in order to 
separate $E_i$'s from $C$.
Then $C^2 \le 2$ and by Lemma 2.2, this is impossible.  
\qed

It turns out that in every case $K+B$ has a 7-complement. Moreover,
we can choose $g$ so that in every case $B_1^\prime$ has only 
one component which is 
a fibre of $\F_2$.

\bigskip
\noindent{\bf Table}

Thus we get the following table. Here,
\begin{itemize}
\item The first column shows the configuration on $S^{\text{min}}$ of the exceptional
curve $E_i$'s, $(-1)$-curves, and the components of $B$. `$\circ$' denote
$(-1)$-curve, `$\bullet$' are the $E_i$'s with self intersection number 
($\le -2$) attached, with `$\leftarrow$' indicating (one possible)
$\Sigma \subset\F_2$ after a suitable sequence of 
contractions of $(-1)$-curves. Squares are curves with non negative
self intersection. 

\item The second column gives the fractional part $B_1$ of the boundary $B$,
or rather, of $D$.

\item The third column gives the number
$(\tfrac67 \le)\hspace{0.1in} max\{b|K+bC+B_1\le 0\}\hspace{0.1in} (<1)$ 
\item The fourth column gives an example of $n$-complements. 
\item The last column lists numerical relations between some relevant divisors on S, with $H$ being the generator of $Pic(S)$.
\end{itemize}

Note that we can compute intersection numbers on $S^{\text{min}}$ using
the crepant pullbacks, and a divisor on $S$ is Cartier iff its crepant
pullback is Cartier, i.e. iff it is integral ({\it cf.}[Sakai]). 

The table
is organized according to $S_C$, the image of $\uppermin{S}$ under the
morphism defined by the linear system $|C|$.

\bigskip
(1) $S=S(A-1)=\Q_2$ (=quadratic cone $\subset \PP^3$, $S_C$ = its Veronese image)
\bigskip

\tablehead
\centering{1} &
\includefigure{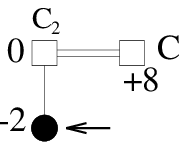} &
\centering{$\frac12C_2$} &
\centering{$\dfrac78$} &
\begin{minipage}{1.4in}$\bullet$ 7-compl.$=0$\\
($K+\frac67C+\frac47C_2\equiv0$)\\
  $\bullet$ trivial 8-compl. \end{minipage}&
\begin{minipage}{1.2in} \begin{center}
$$\begin{array}{ll}
  -K&\equiv C \equiv 2H \\
  C_2&\equiv\frac12H
\end{array} $$ \end{center} \end{minipage}
\\
\hline
\end{tabular}

\bigskip
(2) $S_C=S_7$ (= a del Pezzo with degree 7)
\bigskip

\settablewidth
\centering{2} &
\includefigure{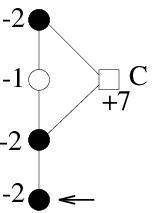} &
\centering{$0$}&
\centering{$\dfrac67$}&
\begin{minipage}{1.4in}\centering{trivial 7-compl.\\
($K+\frac67C\equiv 0$)} \end{minipage} &
\begin{minipage}{1.2in}\begin{center}
$$\begin{array}{ll}
  -K& \equiv H \\
  C & \equiv \frac76H 
\end{array} $$ \end{center} \end {minipage}
\\
\hline
\end{tabular}

\newpage
\par
\bigskip
(3) $S_C = S(A_1+A_2)$
\bigskip

\tablehead

\centering{ 3 }& &
\centering{$\frac12C_2$} &
\centering{$\dfrac9{10}$} &
\begin{minipage}{1.4in}\centering{trivial 10-compl.\\
($K+\frac9{10}C+\frac12C_2\equiv0$)}\end{minipage} &
\begin{minipage}{1.2in} \begin{center}  
$$ \begin{array}{ll}
                 -K  &\equiv C \equiv H\\
                 C_2&\equiv\frac16H          
\end{array} $$ \end{center}
\end{minipage}
\\
\cline{3-5}

& 
\includefigure{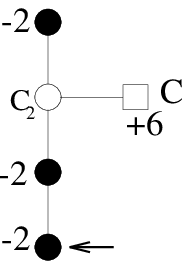} &
\centering{$\frac23C_2$} &
\centering{$\dfrac89$} &
\begin{minipage}{1.4in}\centering{
trivial 9-compl.\\($K+\frac89C+\frac23C_2\equiv0$)} \end{minipage} &

\\
\cline{3-5}
& &
\begin{center}$\frac34C_2$\end{center} &
\begin{center}$\dfrac78$\end{center} &
\begin{center}trivial 8-compl.
\\($K+\frac78C+\frac34C_2\equiv0$)\end{center} &

\\
\cline{3-5}
& &
\begin{center}$\frac45C_2$\end{center} &
\begin{center}$\dfrac{13}{15}$\end{center} &
\begin{center}7-compl.=0\\($K+\frac67C+\frac67C_2\equiv0$)\end{center} &

\\
\cline{3-5}
& &
\begin{center}$\frac56C_2$\end{center} &
\begin{center}$\dfrac{31}{36}$\end{center} &
\begin{center}7-compl.=0\\($K+\frac67C+\frac67C_2\equiv0$)
\end{center} &
\\

\hline

\centering{ 4 }&
\includefigure{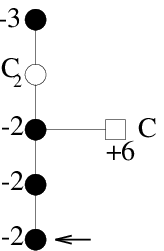} &
\centering{$0$ } &
\centering{$\dfrac89$} &
\begin{minipage}{1.4in}
\centering{$\bullet$ 7-compl.=$C_2$\\($K+\frac67C+\frac27C_2\equiv0$)\\
\bigskip $\bullet$ trivial 9-compl.} \end{minipage}&
\begin{minipage}{1.2in}
\centering{     $ -K  \equiv\frac{8}{12}H$\\
                $  C   \equiv\frac{9}{12}H$\\
                $ C_2 \equiv\frac1{12}H$}
\end{minipage}
\\\hline

\centering{ 5 }&
\includefigure{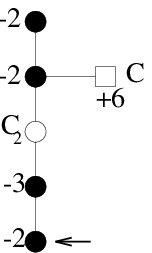} &
\centering{$0$ } &
\centering{$\dfrac9{10}$} &
\begin{minipage}{1.4in}
\centering{$\bullet$ 7-compl.=$3C_2$\\($K+\frac67C+\frac37C_2\equiv0$)\\
\bigskip $\bullet$
trivial 10-compl.}\end{minipage} &
\begin{minipage}{1.2in} \begin{center} 
$$ \begin{array}{ll}
                 -K  &\equiv\frac{9}{15}H\\
                 C  &\equiv\frac{10}{15}H\\
                 C_2&\equiv\frac1{15}H          
\end{array} $$ \end{center}
\end{minipage}
\\\hline

\centering{ 6 }&
\includefigure{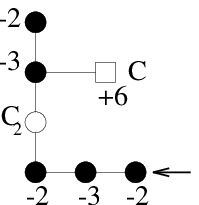} &
\centering{$0$ } &
\centering{$\dfrac78$} &
\begin{minipage}{1.4in}
\centering{$\bullet$ 7-compl.=$2C_2$\\($K+\frac67C+\frac27C_2\equiv0$)\\
\bigskip $\bullet$ trivial 8-compl.}\end{minipage} &
\begin{minipage}{1.2in} \begin{center} 
$$ \begin{array}{ll}
                 -K  &\equiv\frac{14}{40}H\\
                 C  &\equiv\frac{16}{40}H\\
                 C_2&\equiv\frac1{40}H          
\end{array} $$ \end{center}
\end{minipage}
\\\hline
\end{tabular}
\newpage
\tablehead
\centering{ 7 }&
\includefigure{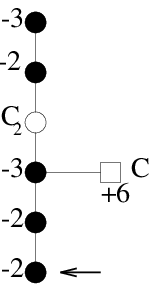} &
\centering{$0$ } &
\centering{$\dfrac{13}{15}$} &
\begin{minipage}{1.4in}\centering{
 7-compl.=$C_2$\\($K+\frac67C+\frac17C_2\equiv0$)
} \end{minipage}&
\begin{minipage}{1.2in} \begin{center} 
$$ \begin{array}{ll}
                 -K  &\equiv\frac{13}{35}H\\
                 C  &\equiv\frac{15}{35}H\\
                 C_2&\equiv\frac1{35}H          
\end{array} $$ \end{center}
\end{minipage}
\\\hline
\centering{ 8 }&
\includefigure{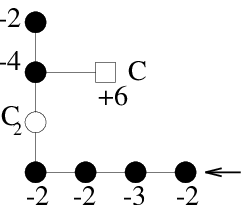} &
\centering{$0$ } &
\centering{$\dfrac{19}{22}$} &
\centering{7-compl.=$C_2$\\($K+\frac67C+\frac17C_2\equiv0$)} &
\begin{minipage}{1.2in} \begin{center} 
$$ \begin{array}{ll}
                 -K  &\equiv\frac{19}{77}H\\
                 C  &\equiv\frac{22}{77}H\\
                 C_2&\equiv\frac1{77}H          
\end{array} $$ \end{center}
\end{minipage}
\\\hline
\end{tabular}
\par
\bigskip
(4) $S_C=S(A_4)$
\bigskip

\begin{tabular}{|p{0.2in}|p{1.2in}|p{0.8in}|p{0.5in}|p{1.4in}|p{1.2in}|} 
\hline
\centering{ 9 }&
\includefigure{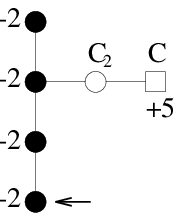} &
\centering{$\frac12C_2$ } &
\centering{$\dfrac{9}{10}$} &
\begin{minipage}{1.4in}
\centering{$\bullet$ 7-compl.=$C_2$\\($K+\frac67C+\frac57C_2\equiv0$)\\
\bigskip $\bullet$ trivial 10-compl.}\end{minipage} &
\begin{minipage}{1.2in} \begin{center} 
$$ \begin{array}{ll}
                 -K  &\equiv C \equiv H\\
                 C_2&\equiv\frac15H          
\end{array} $$ \end{center}
\end{minipage}
\\\cline{3-5}
& &
\begin{center}$\frac23C_2$ \end{center} &
\begin{center}$\dfrac{13}{15}$\end{center} &
\begin{center}7-compl.=0\\($K+\frac67C+\frac57C_2\equiv0$)\end{center}&
\\\hline


\centering{ 10 } &
\includefigure{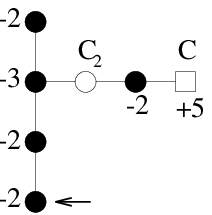} &
\centering{$0$ } &
\centering{$\dfrac{10}{11}$} &
\begin{minipage}{1.4in}
\centering{$\bullet$ 7-compl.=$4C_2$\\($K+\frac67C+\frac47C_2\equiv0$)\\
\bigskip $\bullet$ trivial 11-compl.} \end{minipage}&
\begin{minipage}{1.2in} \begin{center} 
$$ \begin{array}{ll}
                 -K  &\equiv \frac{10}{22}H\\
                 C  &\equiv \frac{11}{22}H\\
                 C_2&\equiv\frac1{22}H          
\end{array} $$ \end{center}
\end{minipage}
\\\cline{3-5}
& &
\begin{center}$\frac12C_2$ \end{center} &
\begin{center}$\dfrac{19}{22}$\end{center} &
\begin{center}7-compl.=0\\($K+\frac67C+\frac47C_2\equiv0$)\end{center}&
\\\hline
\centering{ 11 }&
\includefigure{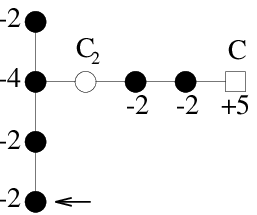} &
\centering{$0$ } &
\centering{$\dfrac{15}{17}$} &
\begin{minipage}{1.4in}
\centering{7-compl.=$3C_2$\\($K+\frac67C+\frac37C_2\equiv0$)}
\end{minipage} &
\begin{minipage}{1.2in}  
$$ \begin{array}{ll}
                 -K  &\equiv\frac{15}{3\cdot 17}H\\
                 C  &\equiv\frac{17}{3\cdot17}H\\
                 C_2&\equiv\frac1{3\cdot17}H          
\end{array} $$
\end{minipage}
\\
\hline
\end{tabular}
\newpage
\tablehead
\centering{12}&
\includefigure{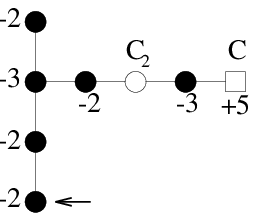} &
\centering{$0$} &
\centering{$\dfrac78$} &
\begin{minipage}{1.4in}
\centering{$\bullet$ 7-compl.=$2C_2$\\($K+\frac67C+\frac27C_2\equiv0$)\\
\bigskip $\bullet$ trivial 8-compl.} \end{minipage}&
\begin{minipage}{1.2in} \begin{center} 
$$ \begin{array}{ll}
                 -K  &\equiv\frac{14}{48}H\\
                 C  &\equiv\frac{16}{48}H\\
                 C_2&\equiv\frac1{48}         
\end{array} $$ \end{center}
\end{minipage}
\\
\hline

\centering {13}&
\includefigure{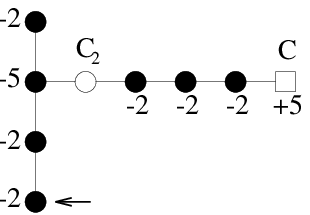} &
\centering{$0$} &
\centering{$\dfrac{20}{23}$}&
\begin{minipage}{1.4in}
\centering{7-compl.=$2C_2$\\($K+\frac67C+\frac27C_2\equiv0$)}
\end{minipage} &
\begin{minipage}{1.2in}  \begin{center}
$$ \begin{array}{ll}
                 -K  &\equiv\frac{20}{4\cdot23}H\\
                 C  &\equiv\frac14H\\
                 C_2&\equiv\frac1{4\cdot23}H         
\end{array} $$ \end{center}
\end{minipage}
\\
\hline

\centering {14}&
\includefigure{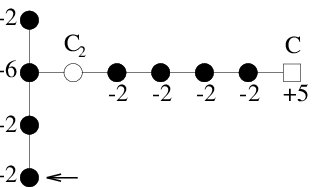} &
\centering{$0$} &
\centering{$\dfrac{25}{29}$}&
\begin{minipage}{1.4in}
\centering{7-compl.=$C_2$\\($K+\frac67C+\frac17C_2\equiv0$)}
\end{minipage} &
\begin{minipage}{1.2in}  \begin{center}
$$ \begin{array}{ll}
                 -K  &\equiv\frac{25}{5\cdot29}H\\
                 C  &\equiv\frac15H\\
                 C_2&\equiv\frac1{5\cdot29}H         
\end{array} $$ \end{center}
\end{minipage}
\\
\hline
\end {tabular}
\par
\smallskip
(5) $S_C=S(D_5)$
\smallskip

\begin{tabular}{|p{0.2in}|p{1.2in}|p{0.8in}|p{0.5in}|p{1.4in}|p{1.2in}|} 
\hline
\centering {15}&
\includefigure{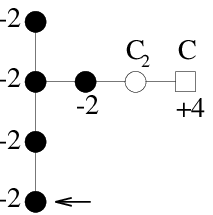}&
\centering{$\frac12 C_2$} &
\centering{$\dfrac78$}&
\begin{minipage}{1.4in}
\centering{$\bullet$ 7-compl.=0\\($K+\frac67C+\frac47C_2\equiv0$)\\
\bigskip $\bullet$ trivial 8-compl.} \end{minipage}&
\begin{minipage}{1.2in}  \begin{center}
$$ \begin{array}{ll}
                 -K  &\equiv C \equiv H\\
                 C_2&\equiv\frac14H         
\end{array} $$ \end{center}
\end{minipage}
\\\hline
\centering {16}&
\includefigure{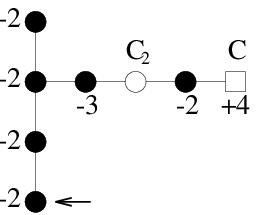} &
\centering{$0$} &
\centering{$\dfrac{8}{9}$}&
\begin{minipage}{1.4in}
\centering{$\bullet$ 7-compl.=$2C_2$\\($K+\frac67C+\frac27C_2\equiv0$)\\
\bigskip $\bullet$ trivial 9-compl.} \end{minipage} &
\begin{minipage}{1.2in}  \begin{center}
$$ \begin{array}{ll}
                 -K  &\equiv \frac8{18}H\\
                  C  &\equiv \frac9{18}H\\
                 C_2&\equiv\frac1{18}H         
\end{array} $$ \end{center}
\end{minipage}
\\\hline
\end {tabular}
\par
\bigskip
(6) $S_C=S(A_3+2A_1)$
\bigskip

\begin{tabular}{|p{0.2in}|p{1.2in}|p{0.8in}|p{0.5in}|p{1.4in}|p{1.2in}|} 
\hline
\centering {17}&
\includefigure{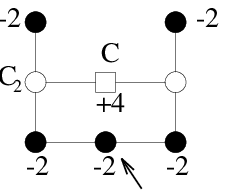} &
\centering{$\frac12 C_2$} &
\centering{$\dfrac{7}{8}$}&
\begin{minipage}{1.4in}
\centering{$\bullet$ 7-compl.=0)\\($K+\frac67C+\frac47C_2\equiv0$)\\
\bigskip $\bullet$ trivial 8-compl.}\end{minipage} &
\begin{minipage}{1.2in}  \begin{center}
$$ \begin{array}{ll}
                 -K  &\equiv C \equiv H\\
                 C_2&\equiv\frac14H         
\end{array} $$ \end{center}
\end{minipage}
\\
\hline

\centering {18}&
\includefigure{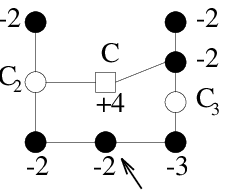} &
\centering{$0$} &
\centering{$\dfrac{6}{7}$}&
\begin{minipage}{1.4in}
\centering{trivial 7-compl.\\($K+\frac67C\equiv0$)}
\end{minipage} &
\begin{minipage}{1.2in}  \begin{center}
$$ \begin{array}{ll}
                 -K  &\equiv \frac{12}{42}H\\
                 C  &\equiv \frac{14}{42}H\\
                 C_2&\equiv\frac3{42}H\\
                 C_3&\equiv\frac2{42}H         
\end{array} $$ \end{center}
\end{minipage}
\\
\hline
\end{tabular}
\newpage
\par\bigskip
(7) $S_C=S(E_6)$
\bigskip

\begin{tabular}{|p{0.2in}|p{1.2in}|p{0.8in}|p{0.5in}|p{1.4in}|p{1.2in}|} 
\hline

\centering {19}&
\includefigure{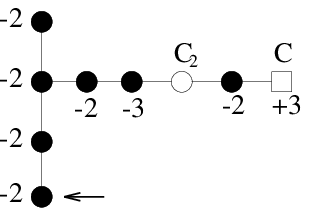}&
\centering{$0$} &
\centering{$\dfrac{6}{7}$}&
\begin{minipage}{1.4in}
\centering{trivial 7-compl.\\($K+\frac67C\equiv0$)} \end{minipage}&
\begin{minipage}{1.2in}  \begin{center}
$$ \begin{array}{ll}
                 -K  &\equiv\frac6{14}H\\
                 C&\equiv\frac7{14}H \\  
                 C_2&\equiv\frac1{14}H      
\end{array} $$ \end{center}
\end{minipage}
\\
\hline
\end {tabular}
\par
\bigskip
(8) $S_C=S(A_5+A_1)$
\bigskip

\begin{tabular}{|p{0.2in}|p{1.2in}|p{0.8in}|p{0.5in}|p{1.4in}|p{1.2in}|} 
\hline
\centering {20}&
\includefigure{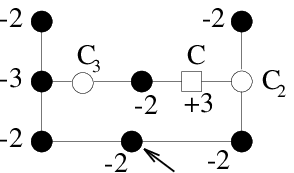}&
\centering{$0$} &
\centering{$\dfrac{6}{7}$}&
\begin{minipage}{1.4in}
\centering{trivial 7-compl.\\($K+\frac67C\equiv0$)} \end{minipage}&
\begin{minipage}{1.2in}  \begin{center}
$$ \begin{array}{ll}
                 -K  &\equiv\frac6{14}H\\
                 C&\equiv\frac7{14}H \\
                 C_2&\equiv \frac2{14}H\\
                 C_3&\equiv \frac1{14}H        
\end{array} $$ \end{center}
\end{minipage}
\\
\hline
\end{tabular}
\newline
\newline

\end{document}